\documentclass{article}
\usepackage{amsfonts}
\usepackage{amsmath}
\usepackage{graphics}
\usepackage{graphicx}
\usepackage{epstopdf}

\input{tcilatex}
\begin{document}

\begin{center}
\bigskip

\textbf{MAGNETIC CURVES ASSOCIATED TO KILLING VECTOR FIELDS IN A GALILEAN
SPACE}

\bigskip

\textbf{Muhittin Evren AYDIN}

Department of Mathematics, Faculty of Science, Firat University, Elazig,
23119, Turkey, meaydin@firat.edu.tr
\end{center}

\textbf{Abstract. }In this paper, we completely classify the magnetic curves
(also $N-$magnetic curves with constant curvature) in a Galilean 3-space
associated to a Killing vector field\linebreak\ $V=v_{1}\partial
_{x}+v_{2}\partial _{y}+v_{3}\partial _{z},$ with $v_{1},v_{2},v_{3}\in 
\mathbb{R}
.$

\bigskip

\textbf{Key words: }Magnetic curve, Lorentz equation, Killing vector field,
cylindrical helix, Galilean space.

\textbf{AMS Subject Classification: }53A35, 53B25, 53Z05.

\section{\textbf{Introduction}}

A closed 2-form $F$ on a complete Riemannian manifold $\left( M,g\right) $
is called a \textit{magnetic field}. The \textit{Lorentz force} of a
magnetic background $\left( M,g,F\right) $ is the skew symmetric $\left(
1,1\right) -$type tensor field $\Phi $ on $M$ satisfying 
\begin{equation}
g\left( \Phi \left( X\right) ,Y\right) =F\left( X,Y\right) ,  \tag{1.1}
\end{equation}%
for any $X,Y\in \Gamma \left( TM\right) .$ Thus a \textit{trajectory}
described by a charged particle moving in a magnetic field $F$ (or a \textit{%
magnetic curve} associated to $F$) is a smooth curve $\gamma $ on $M$
satisfying the \textit{Lorentz equation} (\textit{Newton equation }or\textit{%
\ Landau--Hall equation}):%
\begin{equation}
\nabla _{\dot{\gamma}}\dot{\gamma}=\Phi \left( \dot{\gamma}\right) , 
\tag{1.2}
\end{equation}%
where $\nabla $ is the Levi-Civita connection of $g.$ Note that the Lorentz
force is divergence free, $\func{div}\Phi =0.$

The fact that the Lorentz force is skew symmetric yields a basic property of
magnetic curves, i.e. the following conservation law: particles evolve with
constant speed (and so with constant energy) along the magnetic
trajectories. If the magnetic curve $\gamma $ has unit speed, then it is
called a \textit{normal magnetic curve}. In the sequel, and all over this
paper we study only the unit speed curves.

The study of magnetic fields and their trajectories on Riemannian
(semi-Riemannian) manifolds, situated at the interplay between physics and
differential geometry, has great interest. First problem regarding this
phenomenal field was treated on Riemannian surfaces (see e.g. \cite{8,27}),
then in 3-dimensional context, on $\mathbb{E}^{3}$ \cite{16}$,$ $\mathbb{E}%
_{1}^{3}$ \cite{17}$,$ $\mathbb{S}^{3}$ \cite{9}$,$ $\mathbb{S}^{2}\times 
\mathbb{R}$ \cite{23} etc. For the study of the magnetic curves associated
to magnetic fields on arbitrary dimensional K\"{a}hler, contact manifolds
and Walker 3-manifolds, we also refer the reader to \cite%
{1,2,7,10,12,13,17,18,21,22,25}.

Our aim is to investigate magnetic curves in the Galilean space $\mathbb{G}%
_{3}$ which is one model of the real Cayley-Klein geometries. We classify
the magnetic curves (also $N-$magnetic curves with constant curvature)
associated to the Killing vector field $V=v_{1}\partial _{x}+v_{2}\partial
_{y}+v_{3}\partial _{z}$ in the Galilean 3-space $\mathbb{G}_{3}.$

\section{\textbf{Preliminaries}}

The Galilean geometry is one model of the real Cayley-Klein geometries (for
details see \cite{26}), which has projective signature $(0,0,+,+)$. The 
\textit{absolute figure} of the Galilean geometry is an ordered triple $%
\left\{ \omega ,f,I\right\} $, where $\omega $ is the ideal (absolute)
plane, $f$ a line in $\omega $ and $I$ is the fixed eliptic involution of
the points of $f$. Detailed properties of Galilean space may be found in 
\cite{3}-\cite{6}, \cite{15,19,20,24}.

A plane is called \textit{Euclidean} if it contains $f$, otherwise it is
called\textit{\ isotropic}, i.e., planes $x$ $=const$. are Euclidean, in
particular the plane $\omega $. Other planes are isotropic. A vector $%
X=\left( x_{1},x_{2},x_{3}\right) $ is said to be \textit{non-isotropic }%
(resp. \textit{isotropic}) if $x_{1}\neq 0$ (resp. $x_{1}=0$).

The \textit{Galilean scalar product }between two vectors $X=\left(
x_{1},x_{2},x_{3}\right) $ and $Y=\left( y_{1},y_{2},y_{3}\right) $ is given
by%
\begin{equation*}
\left\langle X,Y\right\rangle _{\mathbb{G}}=\left\{ 
\begin{array}{l}
x_{1}y_{1},\text{ \ \ \ \ \ \ \ \ \ when }x_{1}\neq 0\text{ or }y_{1}\neq 0
\\ 
x_{2}y_{2}+x_{3}y_{3},\text{ when }x_{1}=y_{1}=0.%
\end{array}%
\right.
\end{equation*}%
The norm of the non-isotropic vector $X$ is defined by $\left\Vert
X\right\Vert _{\mathbb{G}}=\left\vert x_{1}\right\vert $; in case when $X$
is isotropic, $\left\Vert X\right\Vert _{\mathbb{G}}=\sqrt{%
x_{2}^{2}+x_{3}^{2}}$. $X$ is called a \textit{unit vector} if $\left\Vert
X\right\Vert _{\mathbb{G}}=1.$

The \textit{cross product} in the sense of Galilean space is (see \cite{3}) 
\begin{equation}
X\times _{\mathbb{G}}Y=\left\{ 
\begin{array}{l}
\left( 0,-%
\begin{vmatrix}
x_{1} & x_{3} \\ 
y_{1} & y_{3}%
\end{vmatrix}%
,%
\begin{vmatrix}
x_{1} & x_{2} \\ 
y_{1} & y_{2}%
\end{vmatrix}%
\right) ,\text{ when }x_{1}\neq 0\text{ or }y_{1}\neq 0 \\ 
\left( 
\begin{vmatrix}
x_{2} & x_{3} \\ 
y_{2} & y_{3}%
\end{vmatrix}%
,0,0\right) ,\text{ \ \ \ \ \ \ \ \ \ \ \ \ \ when }x_{1}=y_{1}=0.%
\end{array}%
\right.  \tag{2.1}
\end{equation}

A smooth admissible curve $\gamma $ (without isotropic tangent vectors) can
be parametrized in $\mathbb{G}_{3}$ by%
\begin{equation}
\gamma :I\subset 
\mathbb{R}
\longrightarrow \mathbb{G}_{3},\text{ }s\longmapsto \left( s,y\left(
s\right) ,z\left( s\right) \right) ,  \tag{2.2}
\end{equation}%
where $s$ is the Galilean invariant of the arc length on $\gamma \left(
I\right) $.

For an admissible curve, the \textit{curvature} $\kappa \left( s\right) $ and%
\textit{\ torsion }$\tau \left( s\right) $ are respectively given by%
\begin{equation*}
\kappa \left( s\right) =\left\Vert \ddot{\gamma}\right\Vert _{\mathbb{G}}=%
\sqrt{\ddot{y}\left( s\right) ^{2}+\ddot{z}\left( s\right) ^{2}}
\end{equation*}%
and%
\begin{equation*}
\tau \left( s\right) =\frac{\det \left( \dot{\gamma},\ddot{\gamma},\dddot{%
\gamma}\right) }{\kappa \left( s\right) ^{2}}.
\end{equation*}%
Obviously, $\kappa \left( s\right) \neq 0,$ $\forall s\in I.$ Thus, the
orthonormal trihedron $\left\{ T,N,B\right\} $ in the sense of Galilean
geometry becomes%
\begin{eqnarray*}
T\left( s\right) &=&\dot{\gamma}\left( s\right) =\left( 1,\dot{y}\left(
s\right) ,\dot{z}\left( s\right) \right) , \\
N\left( s\right) &=&\frac{1}{\kappa \left( s\right) }\left( 0,\ddot{y}\left(
s\right) ,\ddot{z}\left( s\right) \right) , \\
B\left( s\right) &=&\frac{1}{\kappa \left( s\right) }\left( 0,-\ddot{z}%
\left( s\right) ,\ddot{y}\left( s\right) \right) .
\end{eqnarray*}%
In the sequel, the Frenet equations may be expressed by 
\begin{equation*}
\frac{d}{ds}%
\begin{bmatrix}
T \\ 
N \\ 
B%
\end{bmatrix}%
=%
\begin{bmatrix}
0 & \kappa & 0 \\ 
0 & 0 & \tau \\ 
0 & -\tau & 0%
\end{bmatrix}%
\begin{bmatrix}
T \\ 
N \\ 
B%
\end{bmatrix}%
.
\end{equation*}

\section{\textbf{Killing Magnetic Trajectories in }$\mathbb{G}_{3}$}

A smooth vector field $V$ on a Riemannian manifold $\left( M,g\right) $ is a 
\textit{Killing vector field} if the Lie derivative with respect to $V$ of
the metric $g$ vanishes, i.e%
\begin{equation*}
\mathcal{L}_{V}g=0.
\end{equation*}%
A smooth vector field $V$ on $M$ is Killing if and only if it fulfills the
Killing equation%
\begin{equation*}
g\left( \nabla _{X}V,Y\right) +g\left( \nabla _{Y}V,X\right) =0,
\end{equation*}%
for all $X,Y\in \Gamma \left( TM\right) $, where $\nabla$ is the Riemannian
connection on $M$.

The 2-forms on 3-dimensional manifolds may be identified with the
corresponding vector fields via the Hodge $\star $ operator and the volume
form $dv_{g}$ of the manifold. Hence, the magnetic fields correspond to
divergence-free vector fields; in particular, the Killing vector fields
yield an important class of the so called\textit{\ Killing magnetic fields}.

Note that, the cross product of two vector fields $X,Y$ on $M$ can be
defined as%
\begin{equation*}
g\left( X\times Y,Z\right) =dv_{g}\left( X,Y,Z\right) ,
\end{equation*}%
where $X,Y,Z\in \Gamma \left( TM\right) .$ Let $F_{V}=\iota _{V}dv_{g}$ be
the Killing magnetic field corresponding to the Killing vector field $V,$
where $\iota $ denotes the inner product. Equivalently, $F_{V}\wedge $ $%
^{\flat }V=dv_{g},$ where $\flat :TM\longrightarrow TM^{\ast }$ is the
musical isomorphism. Then, the Lorentz force of $F_{V}$ is (see e.g. \cite%
{9,16})%
\begin{equation}
\Phi \left( X\right) =V\times X.  \tag{3.1}
\end{equation}%
Consequently, the relations $\left( 1.2\right) $ and $\left( 3.1\right) $
lead to the Lorentz force of the magnetic background $\left( \mathbb{G}%
_{3},\left\langle ,\right\rangle _{\mathbb{G}},F_{V}\right) :$%
\begin{equation}
\ddot{\gamma}=V\times _{\mathbb{G}}\dot{\gamma},  \tag{3.2}
\end{equation}%
where $V$ is a Killing vector field on $\mathbb{G}_{3}$ and $F_{V}$ the
corresponding Killing magnetic field.

Let $\gamma $ be a curve in $\mathbb{G}_{3},$ parametrized by the arc length
given in the coordinate form%
\begin{equation}
\gamma \left( s\right) =\left( s,y\left( s\right) ,z\left( s\right) \right) ,%
\text{ }s\in I\subset 
\mathbb{R}
,  \tag{3.3}
\end{equation}%
where $y$ and $z$ are smooth functions satisfying the initial conditions:%
\begin{equation}
y\left( 0\right) =y_{0},\text{ }\dot{y}\left( 0\right) =Y_{0}\text{ and\ }%
z\left( 0\right) =z_{0},\text{ }\dot{z}\left( 0\right) =Z_{0}.  \tag{3.4}
\end{equation}

We give a classification of the normal magnetic trajectories associated to
the Killing vector $V=v_{1}\partial _{x}+v_{2}\partial _{y}+v_{3}\partial
_{z}$ in $\mathbb{G}_{3},$ $v_{1},v_{2},v_{3}\in 
\mathbb{R}
.$

\bigskip

\textbf{Theorem 3.1. }\textit{Let }$\gamma $\textit{\ be a normal magnetic
trajectory associated to the Killing vector }$V=v_{1}\partial
_{x}+v_{2}\partial _{y}+v_{3}\partial _{z}$\textit{\ in }$\mathbb{G}_{3}$%
\textit{, satisfying the initial conditions }$\left( 3.4\right) .$\textit{\
Then }$\gamma $\textit{\ has one of the following forms:}

(i) \textit{if }$V$ is isotropic,%
\begin{equation}
\gamma \left( s\right) =\left( s,\frac{v_{3}}{2}s^{2}+Y_{0}s+y_{0},-\frac{%
v_{2}}{2}s^{2}+Z_{0}s+z_{0}\right) ;  \tag{3.5}
\end{equation}

(ii) \textit{otherwise, the cylindrical helix on }$S_{\mathbb{G}}^{1}\left(
r\right) \times l,$\textit{\ where }$S_{\mathbb{G}}^{1}\left( r\right) $ is
a Euclidean circle in $\mathbb{G}_{3}$ with radius\textit{\ }$r=\sqrt{\left( 
\frac{Z_{0}}{v_{1}}-\frac{v_{3}}{v_{1}^{2}}\right) ^{2}+\left( \frac{Y_{0}}{%
v_{1}}-\frac{v_{2}}{v_{1}^{2}}\right) ^{2}}$\textit{\ and }$l$\textit{\ is a
straight line given by }$\left( s,\frac{v_{2}}{v_{1}}s+\left( y_{0}-\frac{%
Z_{0}}{v_{1}}+\frac{v_{3}}{v_{1}^{2}}\right) ,\frac{v_{3}}{v_{1}}s\right. $ $%
\left. +\left( z_{0}+\frac{Y_{0}}{v_{1}}-\frac{v_{2}}{v_{1}^{2}}\right)
\right) ,$\textit{\ parametrized by}%
\begin{eqnarray}
\gamma \left( s\right) &=&\left( s,\frac{\left( Z_{0}-\frac{v_{3}}{v_{1}}%
\right) }{v_{1}}\cos \left( v_{1}s\right) +\frac{\left( Y_{0}-\frac{v_{2}}{%
v_{1}}\right) }{v_{1}}\sin \left( v_{1}s\right) +\frac{v_{2}}{v_{1}}s+\left(
y_{0}-\frac{\left( Z_{0}-\frac{v_{3}}{v_{1}}\right) }{v_{1}}\right) ,\right.
\notag \\
&&\left. \frac{\left( Z_{0}-\frac{v_{3}}{v_{1}}\right) }{v_{1}}\sin \left(
v_{1}s\right) -\frac{\left( Y_{0}-\frac{v_{2}}{v_{1}}\right) }{v_{1}}\cos
\left( v_{1}s\right) +\frac{v_{3}}{v_{1}}s+\left( z_{0}+\frac{\left( Y_{0}-%
\frac{v_{2}}{v_{1}}\right) }{v_{1}}\right) \right) . \quad\qquad(3.6)  \notag
\end{eqnarray}

\bigskip

\textbf{Proof. }The normal magnetic trajectories $\gamma $ are the solutions
of the Lorentz equation $\left( 3.2\right) $. We divide the proof in two
cases:

\bigskip

\textbf{Case 1.} $V$ is isotropic. Using the\ Galilean cross product $\left(
2.1\right) $, we get%
\begin{equation}
\left\{ 
\begin{array}{c}
\ddot{y}=v_{3}, \\ 
\ddot{z}=-v_{2}.%
\end{array}%
\right.  \tag{3.7}
\end{equation}%
After considering the initial conditions $\left( 3.4\right) $ into $\left(
3.7\right) ,$ the normal magnetic trajectory $\gamma $ takes the form $%
\left( 3.5\right) ,$ which gives the statement (i).

\textbf{Case 2. }$V$ is non-isotropic. Then it follows from $\left(
2.1\right) $ and $\left( 3.2\right) $ that 
\begin{equation}
\left\{ 
\begin{array}{l}
\ddot{y}=v_{3}-v_{1}\dot{z}, \\ 
\ddot{z}=v_{1}\dot{y}-v_{2}.%
\end{array}%
\right.  \tag{3.8}
\end{equation}%
We may formulate the Cauchy problem associated to system $\left( 3.8\right) $
and the initial conditions $\left( 3.4\right) $ as follows:%
\begin{equation}
\left\{ 
\begin{array}{l}
\dddot{y}=-v_{1}^{2}\dot{y}+v_{1}v_{2}, \\ 
\dddot{z}=v_{1}v_{3}-v_{1}^{2}\dot{z}.%
\end{array}%
\right.  \tag{3.9}
\end{equation}%
After solving $\left( 3.9\right) ,$ we derive 
\begin{eqnarray*}
y\left( s\right) &=&\frac{v_{2}}{v_{1}}s+\frac{\left( Z_{0}-\frac{v_{3}}{%
v_{1}}\right) }{v_{1}}\cos \left( v_{1}s\right) +\frac{\left( Y_{0}-\frac{%
v_{2}}{v_{1}}\right) }{v_{1}}\sin \left( v_{1}s\right) +\left( y_{0}-\frac{%
\left( Z_{0}-\frac{v_{3}}{v_{1}}\right) }{v_{1}}\right) , \\
z\left( s\right) &=&\frac{v_{3}}{v_{1}}s+\frac{\left( Z_{0}-\frac{v_{3}}{%
v_{1}}\right) }{v_{1}}\sin \left( v_{1}s\right) -\frac{\left( Y_{0}-\frac{%
v_{2}}{v_{1}}\right) }{v_{1}}\cos \left( v_{1}s\right) +\left( z_{0}+\frac{%
\left( Y_{0}-\frac{v_{2}}{v_{1}}\right) }{v_{1}}\right) .
\end{eqnarray*}%
which constitutes the normal magnetic curves we are looking for. It is
straightforward to prove that $\gamma $ is a cylindrical helix wrapped
around $S_{\mathbb{G}}^{1}\left( r\right) \times l$\textit{.}

\bigskip

\textbf{Example 3.1. }Let $\gamma $ be a normal magnetic trajectory
associated to the Killing vector $V=v_{2}\partial _{y}+v_{3}\partial _{z}$
in $\mathbb{G}_{3}.$ It is then expressed by $\left( 3.5\right) .$ Choosing $%
Y_{0}=5,$ $y_{0}=1,$ $Z_{0}=3,$ $z_{0}=4$ and $s\in I=\left[ 0,\pi \right] ,$
the normal magnetic trajectories $\gamma $ becomes the curves in blue color
for $v_{2}=v_{3}=0$, in green color for $v_{2}=v_{3}=1$ and in red color for 
$v_{2}=v_{3}=2$ as in Fig. 1:%
\begin{gather*}
\FRAME{itbpF}{1.9787in}{1.0974in}{0in}{}{}{Figure}{\special{language
"Scientific Word";type "GRAPHIC";maintain-aspect-ratio TRUE;display
"USEDEF";valid_file "T";width 1.9787in;height 1.0974in;depth
0in;original-width 1.9398in;original-height 1.0637in;cropleft "0";croptop
"1";cropright "1";cropbottom "0";tempfilename
'NYMVAW00.wmf';tempfile-properties "XPR";}} \\
\left. 
\begin{array}{l}
\text{\textbf{Fig. 1.}\textit{\ The magnetic trajectories associated to the
Killing vectors} } \\ 
V\in \left\{ (0,0,0),\left( 0,1,1\right) ,\left( 0,2,2\right) \right\} .%
\end{array}%
\right.
\end{gather*}%
\linebreak

\section{$N-$\textbf{Magnetic (}$B-$\textbf{Magnetic) Curves in }$\mathbb{G}%
_{3}$}

In \cite{11}$,$ Bozkurt et al. introduced a new kind of magnetic curves
called $N-$magnetic curves ($B-$magnetic curves) in oriented 3-dimensional
Riemannian manifolds $\left( M,g\right) $ defined as follows:

\bigskip

\textbf{Definition 4.1.} Let $\gamma :I\subset 
\mathbb{R}
\longrightarrow M$ be a curve in an oriented 3-dimensional Riemannian
manifold $\left( M,g\right) $ and $F$ be a magnetic field on $M$. The curve $%
\alpha $ is an $N-$\textit{magnetic curve }(respectively $B-$\textit{%
magnetic curve}) if the normal vector field $N$ (respectively the binormal
vector field $B$) of the curve satisfies the Lorentz force equation, i.e., $%
\nabla _{\dot{\gamma}}N=\Phi \left( N\right) =V\times N$ (respectively $%
\nabla _{\dot{\gamma}}B=\Phi \left( B\right) =V\times B$).

Several charactarizations of $N-$magnetic curves (of $B-$magnetic curves as
well) on $\left( M,g\right) $ were obtained in terms of the curvatures of
the magnetic curve and the curvature of $\left( M,g\right) $ in \cite{11}.

In this section in order to obtain certain results, we assume that the curve 
$\gamma $ has nonzero constant curvature $\kappa _{0}$.

Let consider the curve $\gamma $ in $\mathbb{G}_{3},$ parametrized by%
\begin{equation*}
\gamma \left( s\right) =\left( s,y\left( s\right) ,z\left( s\right) \right) ,
\end{equation*}%
where $y$ and $z$ are smooth functions satisfying the initial conditions:%
\begin{equation}
y\left( 0\right) =y_{0},\text{ }\dot{y}\left( 0\right) =Y_{0},\text{ }\ddot{y%
}\left( 0\right) =T_{0}\text{ and\ }z\left( 0\right) =z_{0},\text{ }\dot{z}%
\left( 0\right) =Z_{0},\text{ }\ddot{z}\left( 0\right) =U_{0}.  \tag{4.1}
\end{equation}

We classify the $N-$magnetic curves with constant curvature $\kappa \left(
s\right) =\kappa _{0}\neq 0,$ corresponding to the Killing vector $%
V=v_{1}\partial _{x}+v_{2}\partial _{y}+v_{3}\partial _{z}$ in $\mathbb{G}%
_{3}$, $v_{1},v_{2},v_{3}\in 
\mathbb{R}
.$

\bigskip

\textbf{Theorem 4.1. }\textit{Let }$\gamma $\textit{\ be a normal }$N-$%
\textit{magnetic trajectory with constant curvature }$\kappa _{0}$\textit{\
associated to the Killing vector }$V=v_{1}\partial _{x}+v_{2}\partial
_{y}+v_{3}\partial _{z}$\textit{\ in }$\mathbb{G}_{3}$\textit{, satisfying
the initial conditions }$\left( 4.1\right) .$\textit{\ Then }$\gamma $%
\textit{\ is one of the following forms:}

(i)\textit{\ If }$v_{1}=v_{2}=v_{3}=0$ ($V=0,$ trivial Killing vector field),%
\textit{\ }%
\begin{equation*}
\gamma \left( s\right) =\left( s,\frac{T_{0}}{2}s^{2}+Y_{0}s+y_{0},\frac{%
U_{0}}{2}s^{2}+Z_{0}s+z_{0}\right) ;
\end{equation*}

(ii) \textit{If }$v_{1}=v_{2}=0$ and $v_{3}\neq 0$\textit{\ }%
\begin{equation*}
\gamma \left( s\right) =\left( s,Y_{0}s+y_{0},\frac{U_{0}}{2}%
s^{2}+Z_{0}s+z_{0}\right) ;
\end{equation*}

(iii)\textit{\ If }$v_{1}=v_{3}=0$ and $v_{2}\neq 0,$ 
\begin{equation*}
\gamma \left( s\right) =\left( s,\frac{T_{0}}{2}%
s^{2}+Y_{0}s+y_{0},Z_{0}s+z_{0}\right)
\end{equation*}

(iv) \textit{If }$v_{1}=0,$ $v_{2}\neq 0,$ $v_{3}\neq 0,$%
\begin{equation*}
\gamma \left( s\right) =\left( s,\frac{U_{0}}{T_{0}}v_{2}s^{2}+Y_{0}s+y_{0},%
\frac{U_{0}}{T_{0}}v_{3}s^{2}+Z_{0}s+z_{0}\right) ;
\end{equation*}

(v) \textit{If }$V$ is \textit{non-isotropic}$,$ $\gamma $ is \textit{the
cylindrical helix on }$S_{\mathbb{G}}^{1}\left( r\right) \times l,$\textit{\
where }$S_{\mathbb{G}}^{1}\left( r\right) $ is a Euclidean circle in $%
\mathbb{G}_{3}$ with radius\textit{\ }$r=\sqrt{\left( \frac{T_{0}}{v_{1}^{2}}%
\right) ^{2}+\left( \frac{U_{0}}{v_{1}^{2}}\right) ^{2}}$\textit{\ and }$l$%
\textit{\ is a straight line given by }$\left( s,\left( Y_{0}-\frac{U_{0}}{%
v_{1}}\right) s+\left( y_{0}+\frac{T_{0}}{v_{1}^{2}}\right) ,\right. $ $%
\left. \left( Z_{0}+\frac{T_{0}}{v_{1}}\right) s+\left( z_{0}+\frac{U_{0}}{%
v_{1}^{2}}\right) \right) $\textit{\ parametrized by}%
\begin{eqnarray*}
\gamma \left( s\right) &=&\left( s,\frac{-T_{0}}{v_{1}^{2}}\cos \left(
v_{1}s\right) +\frac{U_{0}}{v_{1}^{2}}\sin \left( v_{1}s\right) +\left(
Y_{0}-\frac{U_{0}}{v_{1}}\right) s+\left( y_{0}+\frac{T_{0}}{v_{1}^{2}}%
\right) ,\right. \\
&&\left. \frac{-U_{0}}{v_{1}^{2}}\cos \left( v_{1}s\right) -\frac{T_{0}}{%
v_{1}^{2}}\sin \left( v_{1}s\right) +\left( Z_{0}+\frac{T_{0}}{v_{1}}\right)
s+\left( z_{0}+\frac{U_{0}}{v_{1}^{2}}\right) \right) .
\end{eqnarray*}

\bigskip

\textbf{Proof. }Assume that $\gamma $ is a normal $N-$magnetic curve
corresponding to the Killing vector $V=v_{1}\partial _{x}+v_{2}\partial
_{y}+v_{3}\partial _{z}$ in $\mathbb{G}_{3}.$ By the Definition 4.1, we have 
\begin{equation}
\nabla _{\dot{\gamma}}N=\Phi \left( N\right) =V\times _{\mathbb{G}}N, 
\tag{4.2}
\end{equation}%
where $N=\frac{1}{\kappa _{0}}\left( 0,\ddot{y}\left( s\right) ,\ddot{z}%
\left( s\right) \right) .$ If $v_{1}=0,$ it follows from $\left( 4.2\right) $
that%
\begin{equation}
\left\{ 
\begin{array}{l}
\dddot{y}=0, \\ 
\dddot{z}=0, \\ 
v_{2}\ddot{z}-v_{3}\ddot{y}=0.%
\end{array}%
\right.  \tag{4.3}
\end{equation}%
By considering the initial conditions into $\left( 4.1\right) ,$ we
immediately derive the statements (i)-(iv).

In the case $v_{1}\neq 0,$ by using the cross product in $\mathbb{G}_{3}$
and $\left( 4.2\right) ,$ we get that%
\begin{equation}
\left\{ 
\begin{array}{l}
\dddot{y}=-v_{1}\ddot{z}, \\ 
\dddot{z}=v_{1}\ddot{y}.%
\end{array}%
\right.  \tag{4.4}
\end{equation}%
By solving the Cauchy problem associated to system $\left( 4.4\right) $ and
the initial conditions $\left( 4.1\right) ,$ we obtain%
\begin{eqnarray*}
y\left( s\right) &=&\frac{-T_{0}}{v_{1}^{2}}\cos \left( v_{1}s\right) +\frac{%
U_{0}}{v_{1}^{2}}\sin \left( v_{1}s\right) +\left( Y_{0}-\frac{U_{0}}{v_{1}}%
\right) s+\left( y_{0}+\frac{T_{0}}{v_{1}^{2}}\right) , \\
z\left( s\right) &=&\frac{-U_{0}}{v_{1}^{2}}\cos \left( v_{1}s\right) -\frac{%
T_{0}}{v_{1}^{2}}\sin \left( v_{1}s\right) +\left( Z_{0}+\frac{T_{0}}{v_{1}}%
\right) s+\left( z_{0}+\frac{U_{0}}{v_{1}^{2}}\right) ,
\end{eqnarray*}%
which completes the proof.

\bigskip

\textbf{Remark 4.1. }For a trivial magnetic field, $V=0$, i.e. in the case
of vanishing Lorentz force, the magnetic curves are given by the
trajectories of the charged particles moving freely, only under the
influence of gravity. So the solutions of the Lorentz equation become the
geodesics, which satisfy $\nabla _{\dot{\gamma}}\dot{\gamma}=0$. This fact
is not valid for $N-$magnetic curves in $\mathbb{G}_{3}$.

\bigskip

\textbf{Example 4.1. }Let $\gamma $ be a normal $N-$magnetic trajectory
associated to the trivial Killing vector $V=v_{2}\partial _{y}+v_{3}\partial
_{z}$\textit{\ }in $\mathbb{G}_{3}.$ From the first statement of Theorem
4.1, choosing $T_{0}=1,$ $Y_{0}=3,$ $y_{0}=4$ and $U_{0}=1,$ $Z_{0}=2$, $%
z_{0}=1,$ $I=\left[ 0,5\right] ,$ $\gamma \left( s\right) =\left( s,\frac{1}{%
2}u^{2}+3u+4,u^{2}+2u+1\right) .$ Hence the pictures of the $N-$magnetic
trajectories are the curves in blue color, in green color, in red color and
in black color when the Killing vector $V$ is equal to zero, $\partial z,$ $%
\partial y$ and $\partial y+2\partial z,$ respectively, as in Fig. 2:%
\begin{gather*}
\FRAME{itbpF}{1.9061in}{1.9571in}{0in}{}{}{Figure}{\special{language
"Scientific Word";type "GRAPHIC";maintain-aspect-ratio TRUE;display
"USEDEF";valid_file "T";width 1.9061in;height 1.9571in;depth
0in;original-width 1.868in;original-height 1.919in;cropleft "0";croptop
"1";cropright "1";cropbottom "0";tempfilename
'NYMVEX03.wmf';tempfile-properties "XPR";}} \\
\left. 
\begin{array}{l}
\text{\textbf{Fig. 2.}\textit{\ The magnetic trajectories associated to the
Killing vectors} } \\ 
V\in \left\{ (0,0,0),\left( 0,1,0\right) ,\left( 0,0,1\right) ,\left(
0,1,2\right) \right\} .%
\end{array}%
\right.
\end{gather*}

\bigskip

\textbf{Remark 4.2. }Similar results can be obtained by considering $B-$%
magnetic curves (instead of $N-$magnetic curves).

\end{document}